\documentclass{article}
\usepackage[utf8]{inputenc}
\usepackage[letterpaper, portrait, margin=1in]{geometry}
\usepackage{amsmath}
\usepackage{amssymb}
\usepackage{leftidx}
\usepackage{graphicx}
\usepackage{hyperref}
\usepackage{tikz-qtree}
\usepackage{forest}
\usepackage{color}
\usepackage{xcolor}
\usepackage{alltt}

\title{Tor: a toolbox for the continuation of two-dimensional tori in autonomous systems and non-autonomous systems with periodic forcing}
\author{Mingwu Li~\footnote{Email: \href{mailto:mingwli@ethz.ch}{mingwli@ethz.ch}}}
\date{%
    Institute for Mechanical Systems, ETH Z\"{u}rich, Switzerland\\%
}

\begin{document}

\maketitle

\begin{abstract}
We present a toolbox for the continuation of two-dimensional tori in autonomous dynamical systems and non-autonomous systems with periodic forcing. A torus is solved as the solution to a boundary-value problem. Truncated Fourier series and collocation methods are used to discretize the boundary-value problem. This toolbox is based on the continuation package~\href{https://sourceforge.net/projects/cocotools/}{\textsc{coco}}. A few examples are included to demonstrate the functionality of this toolbox.
\end{abstract}

\section{Introduction}
Consider autonomous dynamical system
\begin{equation}
    \frac{dy}{dt}=f(y,p)
\end{equation}
where $f:\mathbb{R}^n\times\mathbb{R}^q\to\mathbb{R}^n$ and non-autonomous system
\begin{equation}
    \frac{dy}{dt}=f(t,y,p)
\end{equation}
where $f:\mathbb{R}\times\mathbb{R}^n\times\mathbb{R}^q\to\mathbb{R}^n$ is periodic in $t$, namely, there exists some positive $T$ such that $f(t,y,p)=f(t+T,y,p)$ for all $t,y$ and $p$. The goal of this report is presenting a unified toolbox \texttt{tor}\footnote{\texttt{tor} toolbox is available at \url{https://github.com/mingwu-li/torus_collocation}} for the continuation of \emph{two-dimensional tori} in the above two types of dynamical systems.

This toolbox includes several core functions as follows
\begin{itemize}
    \item \texttt{ode\_TR2tor}: continuation of tori from a \textit{Neimark-Sacker/torus} bifurcation periodic orbit,
    \item \texttt{ode\_isol2tor}: continuation of tori from an initial solution guess,
    \item \texttt{ode\_tor2tor}: continuation of tori from a saved torus solution,
    \item \texttt{ode\_BP2tor}: continuation of tori along a secondary branch passing through a branch point.
\end{itemize}

This toolbox is based on the continuation package \textsc{coco}~\cite{coco}. There are examples of continuation of tori in \textsc{coco}. Specifically, the continuation of tori in autonomous \emph{Langford} dynamical system is studied in Chapter 9 of~\cite{dankowicz2013recipes}, a book documents the development the package in great detail. In the \texttt{coll} tutorial of \textsc{coco}~\cite{coco-coll}, the continuation of tori in a non-autonomous system with harmonic excitation is studied. Based on these examples, we develop a \emph{general-purpose} toolbox for the continuation of tori in dynamical systems.

In the continuation of periodic orbits, Neimark-Sacker/torus (TR) bifurcation may be detected. It follows that a family of tori emanating from the bifurcation point can be obtained. In \texttt{po} toolbox of \textsc{coco}, one can perform the continuation of periodic orbits. Once a TR bifurcation periodic solution is found using \texttt{po}, one may switch to the continuation of tori starting from such a bifurcation periodic orbit using \texttt{ode\_TR2tor}. Therefore, the \texttt{tor} toolbox connects to the \texttt{po} toolbox seamlessly.

The rest of this report is organized as follows. In section~\ref{sec:prob-form}, problem formulation for tori is given. Specifically, a partial differential equation and phase conditions are derived to solve tori by boundary-value problem approach. Discretization schemes of the boundary-value problem are then discussed in Section~\ref{sec:dist}. In section~\ref{sec:continuation}, parameter continuation of tori is discussed. Following the analysis of dimension deficit of the continuation problem, the initialization of torus solution is presented to construct initial torus solution born from a TR bifurcation periodic orbit. In section~\ref{sec:example}, several examples are presented to demonstrate the effectiveness of the toolbox. Finally, section~\ref{sec:discussion} concludes this report with discussion on the limitation and future enhancement of the current toolbox.

\section{Problem formulation}
\label{sec:prob-form}
\subsection{PDE formulation}
In this section, we focus on the case of non-autonomous dynamical systems. It is straightforward to see that the same formulation holds for autonomous systems as well. The formulation has been derived in~\cite{dankowicz2013recipes,coco-coll}. We include the derivation here for completeness 

For a 2-dimensional quasi-periodic invariant torus, we introduce a torus function $u:\mathbb{S}\times\mathbb{S}\to \mathbb{T}^2$ and two frequencies $\omega_1,\omega_2\in\mathbb{R}$ such that the parallel flow
\begin{equation}
    \frac{d\theta_i}{dt} = \omega_i, \quad i=1,2.
\end{equation}
Substitution of $y(t) = u(\theta_1(t),\theta_2(t))$ into the vector field yields
\begin{equation}
    \omega_1\frac{\partial u}{\partial\theta_1}(\theta_1,\theta_2) + \omega_2\frac{\partial u}{\partial\theta_2}(\theta_1,\theta_2) = f(t,u(\theta_1,\theta_2),p).
\end{equation}
The above equation is a first order \emph{quasi-linear} partial differential equation (PDE). With the method of characteristics, we let $\theta_1 = \varphi + \omega_1 t$ and $\theta_2 = \omega_2t$ and introduce
\begin{equation}
    v(\varphi,t): = u(\varphi + \omega_1 t,\omega_2 t)
\end{equation}
for $\varphi\in\mathbb{S}$ and $t\in [0, {2\pi}/{\omega_2}]$. 
From the transformation $(\varphi,t)\mapsto(\theta_1,\theta_2)$, we obtain $\varphi=\theta_1-\varrho\theta_2$ and $t = {\theta_2}/{\omega_2}$ with $\varrho=\omega_1/\omega_2$. It follows that
\begin{gather}
    \omega_1\frac{\partial u}{\partial\theta_1}(\theta_1,\theta_2) = \omega_1\bigg(\frac{\partial v}{\partial \varphi}\frac{\partial \varphi}{\partial\theta_1} + \frac{\partial v}{\partial t}\frac{\partial t}{\partial\theta_1}\bigg)= \omega_1\frac{\partial v}{\partial \varphi},\\
    \omega_2\frac{\partial u}{\partial\theta_2}(\theta_1,\theta_2) = \omega_2\bigg(\frac{\partial v}{\partial \varphi}\frac{\partial \varphi}{\partial\theta_2} + \frac{\partial v}{\partial t}\frac{\partial t}{\partial\theta_2}\bigg)= -\omega_1\frac{\partial v}{\partial \varphi}+\frac{\partial v}{\partial t},
\end{gather}
and the vector field becomes
\begin{equation}
\label{eq:pde}
    \frac{\partial v}{\partial t} (\varphi,t) =  f(t, v(\varphi,t),p),\quad \varphi\in\mathbb{S},t\in [0,T]
\end{equation}
where $T:={2\pi}/{\omega_2}$.

As for boundary conditions, note that
\begin{gather}
    v(\varphi,0)=u(\varphi,0),\\ u(\varphi+2\pi\varrho,0)=u(\varphi+2\pi\varrho,2\pi)=u(\varphi+\omega_1\cdot2\pi/\omega_2,\omega_2\cdot2\pi/\omega_2)=v(\varphi,2\pi/\omega_2)=v(\varphi,T),
\end{gather}
and then an \emph{all-to-all} coupling condition is obtained as follows
\begin{equation}
\label{eq:pde-bc}
    v(\varphi+2\pi\varrho,0)=v(\varphi,T),\quad \forall\varphi\in\mathbb{S}.
\end{equation}
We have a boundary-value problem (BVP) governed by \eqref{eq:pde} and \eqref{eq:pde-bc}. However, we need phase conditions to yield unique solution to the BVP.

\subsection{Phase conditions}
In the case of non-autonomous systems, we need a phase condition given $\varphi\in\mathbb{S}$. We use Poincar\'e phase condition as follows
\begin{equation}
\label{eq:phase1}
    \langle v_\varphi^\ast(0,0), v(0,0)-v^\ast(0,0)\rangle = 0,
\end{equation}
where $v^\ast(\varphi,t)$ is a \emph{known} function of $\varphi$ and $t$, and $v^\ast_\varphi$ is its partial derivative with respect to $\varphi$. During continuation, $v^\ast(\varphi,t)$ will be updated as the solution of previous continuation step and then the Poincar\'e section is \emph{moving} during continuation. Integral phase condition provides an alternative to Poincar\'e phase condition. The reader may refer to~\cite{schilder2005continuation} for more details about the integral phase condition. 

In the case of autonomous systems, we need one more phase condition given \eqref{eq:pde} becomes autonomous. For consistency, we again use Poincar\'e phase condition
\begin{equation}
\label{eq:phase2}
    \langle v_t^\ast(0,0), v(0,0)-v^\ast(0,0)\rangle = 0,
\end{equation}
where $v^\ast_t$ is the partial derivative of $v^\ast(\varphi,t)$ with respect to $t$. Likewise, the Poincar\'e section here is updated during continuation.

\section{Discretization}
\label{sec:dist}
One need to discretize the domain $\mathbb{S}\times[0,T]$ to obtain approximated BVP. We first apply truncated Fourier series expansion of $v(\varphi,t)$ in $\varphi$, yielding a \emph{multi-segments} BVP with ordinary differential equations (ODEs). Then, collocation method is used to approximate the obtained ODEs to yield a set of nonlinear algebraic equations.

\subsection{Fourier series expansion}
Following~\cite{dankowicz2013recipes}, let $\chi$ represent truncated Fourier series expansion of a component of $v(\varphi,\tau)$, we have
\begin{equation}
\label{eq:FourierExp}
    v_i(\varphi,t)\approx \chi(\varphi,t):=  a_0(t) + \sum_{k=1}^N a_k(t)\cos( k\varphi) + b_k(t)\sin( k\varphi).
\end{equation}
It follows that
\begin{equation}
    \begin{pmatrix}\chi(0,t)\\\chi(\frac{2\pi}{2N+1},t)\\\vdots\\\chi(\frac{4\pi N}{2N+1},t)\end{pmatrix} = \mathcal{F}^{-1} \begin{pmatrix}a_0(t)\\a_1(t)\\b_1(t)\\\vdots\\a_N(t)\\b_N(t)\end{pmatrix}
\end{equation}
where $\mathcal{F}$ is \emph{discrete Fourier transform matrix} and given by equation 9.21 in~\cite{dankowicz2013recipes}.

Likewise, we have
\begin{equation}
    v_i(\varphi+2\pi\varrho,t)\approx \chi(\varphi+2\pi\varrho,t):=  a'_0(t) + \sum_{k=1}^N a'_k(t)\cos( k\varphi) + b'_k(t)\sin( k\varphi).
\end{equation}
It follows that
\begin{equation}
    \begin{pmatrix}a'_0(t)\\a'_1(t)\\b'_1(t)\\\vdots\\a'_N(t)\\b'_N(\tau)\end{pmatrix} = \mathcal{R} \begin{pmatrix}a_0(t)\\a_1(t)\\b_1(t)\\\vdots\\a_N(t)\\b_N(t)\end{pmatrix}
\end{equation}
where $\mathcal{R}$ is given by equation (38) in~\cite{coco-coll} or equation (9.24) in~\cite{dankowicz2013recipes}.

From the derived \emph{all-to-all} coupling boundary condition, we have
\begin{equation}
    \chi(\varphi,T) = \chi(\varphi+2\pi\varrho,0),
\end{equation}
or equivalently
\begin{equation}
    a_0(T) + \sum_{k=1}^N a_k(T)\cos( k\varphi) + b_k(T)\sin( k\varphi) = a'_0(0) + \sum_{k=1}^N a'_k(0)\cos( k\varphi) + b'_k(0)\sin( k\varphi)
\end{equation}
from which we obtain
\begin{equation}
    \begin{pmatrix}a'_0(0)\\a'_1(0)\\b'_1(0)\\\vdots\\a'_N(0)\\b'_N(0)\end{pmatrix} = \begin{pmatrix}a_0(T)\\a_1(T)\\b_1(T)\\\vdots\\a_N(T)\\b_N(T)\end{pmatrix}.
\end{equation}
Substitution yields
\begin{equation}
    \mathcal{R}\begin{pmatrix}a_0(0)\\a_1(0)\\b_1(0)\\\vdots\\a_N(0)\\b_N(0)\end{pmatrix} = \begin{pmatrix}a_0(T)\\a_1(T)\\b_1(T)\\\vdots\\a_N(T)\\b_N(T)\end{pmatrix}\implies \mathcal{R}\mathcal{F}\begin{pmatrix}\chi(0,0)\\\chi(\frac{2\pi}{2N+1},0)\\\vdots\\\chi(\frac{4\pi N}{2N+1},0)\end{pmatrix} = \mathcal{F}\begin{pmatrix}\chi(0,T)\\\chi(\frac{2\pi}{2N+1},T)\\\vdots\\\chi(\frac{4\pi N}{2N+1},T)\end{pmatrix}.
\end{equation}

We proceed to discretize the continuous family $v(\varphi,\tau)$ by restricting attention to the
mesh $\{\varphi_j\}_{j=1}^{2N+1}$, where $\varphi_j:=\frac{2(j-1)\pi}{2N+1}$. Following the above equation, the discrete \emph{all-to-all} boundary condition is given by
\begin{equation}
\label{eq:all-to-all}
    (\mathcal{F}\otimes I_n) \begin{pmatrix}v(\varphi_1,T)\\v(\varphi_2,T)\\\vdots\\v(\varphi_N,T)\end{pmatrix} = ((\mathcal{R}\mathcal{F})\otimes I_n) \begin{pmatrix}v(\varphi_1,0)\\v(\varphi_2,0)\\\vdots\\v(\varphi_N,0)\end{pmatrix},
\end{equation}
where $I_n$ denotes identity matrix of size $n$.

In addition, we obtain the following set of ODEs
\begin{equation}
\label{eq:odes}
    \frac{dv(\varphi_j,t)}{dt} =  f(t,v(\varphi_j,t),p),\quad j=1,\cdots,2N+1, t\in[0,T].
\end{equation}
Along with the boundary condition \eqref{eq:all-to-all} and appropriate discrete phase conditions (which will be discussed in section~\ref{sec:dist-phase}), we have a \emph{multi-segments} BVP with a set of ODEs.

\subsection{Collocation}
We apply collocation methods to solve the above multi-segments BVP. Here we summarize the basic idea of collocation methods. One may refer to~\cite{dankowicz2013recipes} for more details.
\begin{enumerate}
    \item Division of domain: the domain $[0,T]$ is divided into subintervals. 
    \item Approximation of unknown functions: the unknown function within each subinterval is approximated using Lagrangian interpolation of values of the function at a set of base points, and the continuity of the function between two adjacent subintervals is imposed.
    \item Approximation of differential equations: the ODEs are satisfied at a set of collocation nodes.
\end{enumerate}
The collocation method with adaptive mesh has been implemented in \texttt{coll} toolbox of \textsc{coco}. The \texttt{tor} toolbox here is built upon the \texttt{coll} toolbox.

\subsection{Phase conditions}
\label{sec:dist-phase}
The imposition of phase conditions \eqref{eq:phase1} and \eqref{eq:phase2} is straightforward with the above discretization. We have $v^\ast_t(0,0)=f(0,v^\ast(0,0),p)$. To obtain the approximation of $v_\varphi^\ast(0,0)$, let recall \eqref{eq:FourierExp}
and then 
\begin{equation}
    \chi_\varphi(\varphi,t)= \sum_{k=1}^N \left( kb_k(t)\right)\cos(k\varphi)+\left(-ka_k(t)\right)\sin(k\varphi).
\end{equation}
Plugging $(\varphi,t)=(0,0)$ yields 
\begin{equation}
    \chi_\varphi(0,0)=\sum_{k=1}^N kb_k(0).
\end{equation}
Recall
\begin{equation}
     \begin{pmatrix}a_0(0)\\a_1(0)\\b_1(0)\\\vdots\\a_N(0)\\b_N(0)\end{pmatrix} = \mathcal{F}\begin{pmatrix}\chi(\varphi_1,0)\\\chi(\varphi_2,0)\\\vdots\\\chi(\varphi_{2N+1},0)\end{pmatrix},
\end{equation}
we can express $\chi_\varphi(0,0)$ as a linear function of $\{\chi(\varphi_j,0)\}_{j=1}^{2N+1}$. It follows that we can obtain $v_\varphi^\ast(0,0)$ in terms of $\{v^\ast(\varphi_j,0)\}_{j=1}^{2N+1}$.

\section{Parameter continuation}
\label{sec:continuation}
\subsection{Dimension deficit}
\label{sec:deficit}
Recall system parameters $p\in\mathbb{R}^q$, we regard $\omega_1$, $\omega_2$ and $\varrho$ as system parameters as well. Other than the \emph{all-to-all} coupling condition, we have the following boundary conditions
\begin{equation}
    T_0=0, T=\frac{2\pi}{\omega_2},\varrho=\frac{\omega_1}{\omega_2},
\end{equation}
where $T_0$ is the initial time.
With initially inactive system parameters, the dimension deficit is $-1$. We also need to account for phase conditions. It follows that
\begin{itemize}
    \item \emph{autonomous systems}: two phase conditions are included and then the dimension deficit is $-3$. So four parameters need to be released to yield an one-dimensional solution manifold. For instance, one may release $\omega_1,\omega_2$ and $p_{1,2}$ to obtain a family of tori with fixed rotation number $\varrho$,
    \item \emph{non-autonomous systems}: one phase condition is included and we have a coupling condition for system parameters, namely, $\Omega_2-\omega_2=0$, where $\Omega_2$ is a parameter included in $p$ and characterizes the period of non-autonomous forcing ($T=2\pi/\Omega_2$). So the dimension deficit is again $-3$ and we need to release four parameters to obtain a one-dimensional manifold of tori.
\end{itemize}

\subsection{Initialization}
In the continuation of periodic orbits, Neimark-Sacker (torus) bifurcation may be observed. The occurrence of torus bifurcation indicates the born of tori. Here we provide an initial guess to such a torus. The function \texttt{ode\_TR2tor} is based on the initialization presented in this section. The initialization here is adapted from the derivation in~\cite{olikara2010computation}.

Consider dynamical system
\begin{equation}
    \dot{x}=f(t,x,p),\quad x(t_0)=x_0
\end{equation}
with a solution $x(t)$. It follows that
\begin{equation}
    x(t)=x_0+\int_{t_0}^t f(s,x(s),p)\mathrm{d}s
\end{equation}
and hence
\begin{equation}
    \frac{\partial x(t)}{\partial x_0}=\mathbb{I}+\int_{t_0}^t \frac{\partial f(s,x(s),p)}{\partial x}\frac{\partial x(s)}{\partial x_0}\mathrm{d}s.
\end{equation}
Let $\Phi(t,t_0):=D_{x_0}x(t)=\frac{\partial x(t)}{\partial x_0}$, we have
\begin{equation}
    \dot{\Phi}=f_x\Phi,\quad \Phi(t_0)=\mathbb{I}
\end{equation}
and the transition matrix $\mathcal{M}(t,t_0)$ is given by $\Phi(t,t_0)$ and the monodromy matrix is $\mathcal{M}(t_0+T,t_0)$.
In general, we solve the above variational equation with initial condition $\Phi_0$. Let $\delta x_0=\Phi_0 v$, it follows that
\begin{equation}
    \delta x(t)=\Phi(t,t_0) v=\Phi(t,t_0)\Phi_0^{-1}\delta x_0\implies \mathcal{M}(t,t_0)=\Phi(t,t_0)\Phi_0^{-1}.
\end{equation}

Suppose $\mathcal{M}(t_0+T,t_0)$ has an eigenvalue $e^{\mathrm{i}\alpha}$ and an eigenvector $v$ such that 
\begin{equation}
    \mathcal{M}(t_0+T,t_0)v = e^{\mathrm{i}\alpha}v.
\end{equation}
$\mathrm{Re}(v)$ and $\mathrm{Im}(v)$ span an invariant subspace to the monodromy matrix.
Define $u(t)=e^{-\mathrm{i}\alpha (t-t_0)/T}\mathcal{M}(t,t_0)v$, it follows that $u(t_0+T)=u(t_0)=v$ and then $u(t)$ is a periodic function with period $T$. 
For a periodic reference solution $x_p(t)$ with period $T$, we know the coefficient matrix $f_x$ is time-periodic and then the transition matrix $\mathcal{M}(t,t_0)$ satisfies
\begin{equation}
    \mathcal{M}(t+T,t_0)=\mathcal{M}(t+T,t)\mathcal{M}(t,t_0)=\mathcal{M}(t+T,t_0+T)\mathcal{M}(t_0+T,t_0)=\mathcal{M}(t,t_0)\mathcal{M}(t_0+T,t_0)
\end{equation}
where the first two equalities use the semi-group property of state transition matrix and the last one uses the periodicity of coefficient matrix.
Then we can show that $u(t)$ is an eigenvector to $\mathcal{M}(t+T,t)$ with same eigenvalue $e^{\mathrm{i}\alpha}$. Note
\begin{equation}
    \mathcal{M}(t+T,t)\mathcal{M}(t,t_0)v=\mathcal{M}(t,t_0)\mathcal{M}(t_0+T,t_0)v=e^{\mathrm{i}\alpha}\mathcal{M}(t,t_0)v,
\end{equation}
we have $\mathcal{M}(t+T,t)u(t)=e^{\mathrm{i}\alpha}u(t)$.

So we can define a torus function perturbation
\begin{equation}
    \hat{u}(\theta_1,t) = \mathrm{Re}[e^{\mathrm{i}\theta_1}u(t)]=\cos\theta_1\mathrm{Re}(u(t))-\sin\theta_1\mathrm{Im}(u(t)),\quad\theta_1\in\mathbb{S}.
\end{equation}
It follows that
\begin{equation}
    \mathcal{M}(t+T)\hat{u}(\theta_1,t)=\mathrm{Re}[\mathcal{M}(t+T)e^{\mathrm{i}\theta_1}u(t)]=\mathrm{Re}[e^{\mathrm{i}\theta_1}\mathcal{M}(t+T)u(t)]=\mathrm{Re}[e^{\mathrm{i}\theta_1}e^{\mathrm{i}\alpha}u(t)]=\hat{u}(\theta_1+\alpha, t)
\end{equation}
which indicates the rotation of $\hat{u}(\theta_1,t)$ with angle $\alpha$ when state transition is applied. The perturbed torus solution is then given by
\begin{equation}
    \tilde{u}(\theta_1,t) = \mathcal{M}(t,t_0)\hat{u}(\theta_1,t_0) = \mathrm{Re}[e^{\mathrm{i}\theta_1}\mathcal{M}(t,t_0)v] =\hat{u}(\theta_1+\alpha(t-t_0)/T,t).
\end{equation}
So the torus initial solution is given by
\begin{equation}
    \hat{x}(\theta_1,t) = x(t) + \epsilon\tilde{u}(\theta_1,t)=x(t)+\epsilon\hat{u}(\theta_1+\omega_1(t-t_0),t),\quad \forall\theta_1\in\{\varphi_j\}_{j=1}^{2N+1}
\end{equation}
where $\epsilon$ controls the amount of perturbation. We have used the fact that $T=\frac{2\pi}{\omega_2}$ and $\alpha=\omega_1T$.

Note that the eigenvector is still an eigenvector after multiplication by a complex number. So we would like to choose the one such that $\langle v_\mathrm{R},v_\mathrm{I}\rangle=0$, where $v_\mathrm{R}=\mathrm{Re}(v)$ and $v_\mathrm{I}=\mathrm{Im}(v)$. Let the complex number be $e^{\mathrm{i}\theta}$, it follows that
\begin{equation}
    e^{\mathrm{i}\theta}(v_\mathrm{R}+\mathrm{i}v_\mathrm{I})=(\cos\theta+\mathrm{i}\sin\theta)(v_\mathrm{R}+\mathrm{i}v_\mathrm{I})=(\cos\theta v_\mathrm{R}-\sin\theta v_\mathrm{I})+\mathrm{i}(\cos\theta v_\mathrm{I}+\sin\theta v_\mathrm{R}).
\end{equation}
We ask for
\begin{equation}
    \langle\cos\theta v_\mathrm{R}-\sin\theta v_\mathrm{I},\cos\theta v_\mathrm{I}+\sin\theta v_\mathrm{R}\rangle=0,
\end{equation}
from which we obtain
\begin{equation}
    \langle v_\mathrm{R},v_\mathrm{I}\rangle\cos2\theta+0.5(\langle v_\mathrm{R},v_\mathrm{R}\rangle-\langle v_\mathrm{I},v_\mathrm{I}\rangle)\sin2\theta=0\implies \tan2\theta=\frac{2\langle v_\mathrm{R},v_\mathrm{I}\rangle}{\langle v_\mathrm{I},v_\mathrm{I}\rangle-\langle v_\mathrm{R},v_\mathrm{R}\rangle}.
\end{equation}
We can solve $\theta$ from the above equation with \texttt{atan2} function.

\section{Examples}
\label{sec:example}
\subsection{Langford dynamical system}
Consider the \emph{Langford} dynamical system
\begin{gather}
    \dot{x}_1 = (x_3-0.7)x_1-\omega x_2,\\
    \dot{x}_2 = \omega x_1+(x_3-0.7)x_2,\\
    \dot{x}_3 = 0.6+x_3-\frac{1}{3}x_3^3-(x_1^2+x_2^2)(1+\rho x_3)+\epsilon x_3x_1^3.
\end{gather}
Following the analysis in~\cite{dankowicz2013recipes}, the system has a family of periodic orbits when $\epsilon=0$, and such solution undergoes a torus bifurcation at $\rho=\rho^\ast\approx0.615$.

We proceed to use \texttt{po} toolbox in \textsc{coco} to find the family of periodic orbits and detect TR bifurcation point along the solution manifold. Specifically, we first use \texttt{ode45} to perform forward simulation and generate an approximated periodic orbit solution. Then we call \texttt{ode\_isol2po} to encode the boundary-value problem of periodic orbits.
\begin{alltt}
%% construct initial periodic orbit
p0 = [3.5; 1.5; 0]; % [om ro eps]
T  = 2*pi/p0(1);
[~,x0] = ode45(@(t,x) lang(x,p0), 0:100*T, [0.3; 0.4; 0]); % transient
[t0,x0] = ode45(@(t,x) lang(x,p0), linspace(0,T,100), x0(end,:)); % periodic solution
figure;
plot(t0,x0);

%% continuation of periodic orbit
prob = coco_prob();
prob = ode_isol2po(prob, '', @lang, @lang_DFDX, @lang_DFDP, ...
  t0, x0, {'om','rho','eps'}, p0);

coco(prob, 'po', [], 1, 'rho', [0.2 2]);
\end{alltt}
Here \texttt{lang} denotes the vector field of the dynamical system, \texttt{lang\_DFDX} and \texttt{lang\_DFDP} represent the derivative of the vector field with respect to state $x$ and system parameters $p$ respectively. The system parameters in this example are $(\omega,\rho,\epsilon)$ and named as \texttt{om}, \texttt{rho} and \texttt{eps} as above. The vector field and its derivatives are encoded as follows
\begin{alltt}
function y = lang(x, p)

x1  = x(1,:);
x2  = x(2,:);
x3  = x(3,:);
om  = p(1,:);
ro  = p(2,:);
eps = p(3,:);

y(1,:) = (x3-0.7).*x1-om.*x2;
y(2,:) = om.*x1+(x3-0.7).*x2;
y(3,:) = 0.6+x3-x3.^3/3-(x1.^2+x2.^2).*(1+ro.*x3)+eps.*x3.*x1.^3;

end

function J = lang_DFDX(x, p)

x1  = x(1,:);
x2  = x(2,:);
x3  = x(3,:);
om  = p(1,:);
ro  = p(2,:);
eps = p(3,:);

J = zeros(3,3,numel(x1));
J(1,1,:) = (x3-0.7);
J(1,2,:) = -om;
J(1,3,:) = x1;
J(2,1,:) = om;
J(2,2,:) = (x3-0.7);
J(2,3,:) = x2;
J(3,1,:) = -2*x1.*(1+ro.*x3)+3*eps.*x3.*x1.^2;
J(3,2,:) = -2*x2.*(1+ro.*x3);
J(3,3,:) = 1-x3.^2-ro.*(x1.^2+x2.^2)+eps.*x1.^3;

end

function J = lang_DFDP(x, p)

x1 = x(1,:);
x2 = x(2,:);
x3 = x(3,:);

J = zeros(3,size(p,1),numel(x1));
J(1,1,:) = -x2;
J(2,1,:) = x1;
J(3,2,:) = -x3.*(x1.^2+x2.^2);
J(3,3,:) = x3.*x1.^3;

end
\end{alltt}

Indeed, a TR bifurcation periodic orbit is found at $\rho\approx0.61545$. We then move to the continuation of tori in the system. We call \texttt{ode\_TR2tor} to switch from the continuation of periodic orbits to the continuation of tori.
\begin{alltt}
T_po   = 5.3; 
T_ret  = 2*pi/3.5;
varrho = T_ret/T_po;
bd    = coco_bd_read('po');
TRlab = coco_bd_labs(bd, 'TR');
prob = coco_prob();
prob = coco_set(prob, 'cont', 'NAdapt', 5, 'h_min',...
    1e-3, 'PtMX', 50, 'h_max', 10, 'bi_direct', false);
prob = ode_TR2tor(prob, '', 'po', TRlab, 50);

coco(prob, 'tr1', [], 1, {'varrho','rho','om1','om2','eps'},[varrho,0.44]);
\end{alltt}
Here \texttt{varrho}, \texttt{om1} and \texttt{om2} correspond to $\varrho$, $\omega_1$ and $\omega_2$ in the PDE formulation respectively. \texttt{TRlab} gives the label of TR bifurcation solution in \texttt{po} run. The 50 in \texttt{prob = ode\_TR2tor(prob, '', 'po', TRlab, 50)} characterizes the number of segments, namely, $N$ in section~\ref{sec:dist}. Recall the number of segments is given by $2N+1$, so we have 101 segments in this case. By default, $N=10$. Please type \texttt{help ode\_TR2tor} in \textsc{matlab} for more details of the syntax of this constructor.

In the above continuation run \texttt{tr1}, $(\varrho,\rho,\omega_1,\omega_2)$ are varied and a one-dimensional solution manifold is obtained. As discussed in section~\ref{sec:deficit}, four parameters need to be released to yield a one-dimensional manifold. Indeed, \texttt{eps} does not change during the continuation run. As an alternative, we may fix $\varrho$ and then another one-dimensional manifold of tori is obtained. To demonstrate it, we call constructor \texttt{ode\_tor2tor} as follows
\begin{alltt}
%% continuation of torus from previous solution
bd   = coco_bd_read('tr1');
lab  = coco_bd_labs(bd, 'EP');
lab  = max(lab);
prob = coco_prob();
prob = coco_set(prob, 'cont', 'NAdapt', 5, 'h_min',...
    1e-3, 'PtMX', 30, 'h_max', 10, 'bi_direct', true);
prob = ode_tor2tor(prob, '', 'tr1', lab);

coco(prob, 'tr2', [], 1, {'eps','rho','om1','om2','varrho'});
\end{alltt}
Here we start continuation from a solution found in previous run \texttt{tr1}. In current run \texttt{tr2}, $\varrho$ is not varied in continuation, as can be seen in the continuation history.

Once a torus solution is found, you may call \texttt{plot\_torus} to visualize the solution. When $\rho$ is away from $\rho^\ast$, the size of tori born from TR periodic orbit grows, as can be seen in the visualization of a family of tori
\begin{alltt}
%% visualization
figure; coco_plot_bd('tr2','rho','eps');
for lab = 1:5
    plot_torus('','tr1', lab, [1 2 3], 0.75); pause(1);
end
\end{alltt}

To validate an obtained torus solution, one may perform forward simulation with an initial condition on the torus
\begin{alltt}
lab = 5;
plot_torus('','tr1', lab, [1 2 3]);hold on
sol = tor_read_solution('','tr1',lab);
p   = sol.p(1:end-3);
xbp = sol.xbp;
[t,x] = ode45(@(t,x) lang(x,p), 0:0.01:100*T, xbp(1,:,1));
plot3(x(:,1),x(:,2),x(:,3),'r-');
\end{alltt}
Here \texttt{tor\_read\_solution} extracts the result of continuation run. For more details of the arguments of \texttt{plot\_torus} and \texttt{tor\_read\_solution}, please check their corresponding help info.

\subsection{Van der Pol oscillator}
Consider Van der Pol oscillator subject to harmonic excitation
\begin{equation}
    \ddot{x}-c(1-x^2)\dot{x}+x=a\cos\Omega t.
\end{equation}
We are interested in the quasiperiodic response of such a system.

We start by constructing initial solution guess with forward simulation
\begin{alltt}
% pnames = [  om     c     a ]
p0       = [ 1.5111; 0.11; 0.1 ];

T_po = 2*pi; % Approximate period
N    = 10;   % 2N+1 = Number of orbit segments
tout = linspace(0, T_po, 2*N+2);

T_ret = 2*pi/p0(1); % return time = 2*pi/om
tt    = linspace(0,1,10*(2*N+1))';
t1    = T_ret*tt;
x0    = zeros(numel(tt),2,2*N+1);

figure; hold on
for i=1:2*N+1 
% For each point on orbit, flow for return time and reconstitute 3D trajectory
  x1       = [2*cos(tout(i)) 2*sin(tout(i))];
  [~, x1] = ode45(@(t,x) vdp(t,x,p0), 10*t1, x1); % Transient simulation
  [~, x1] = ode45(@(t,x) vdp(t,x,p0), t1, x1(end,:));
  x0(:,:,i) = x1;
  plot3(x1(:,1),t1,x1(:,2),'r-');
end

varrho = T_ret/T_po;
\end{alltt}

We proceed to use constructor \texttt{ode\_isol2tor} to perform continuation of tori with the generated initial solution guess. Given the dynamical system is non-autonomous, \texttt{coco\_set} is called to change default \emph{autonomous} setting. In addition, the excitation frequency parameter $\Omega$ needs to be named as \texttt{Om2} given the coupling condition $\Omega_2-\omega_2$ needs to be imposed, as discussed in section~\ref{sec:deficit}. With these observations, we have
\begin{alltt}
prob = coco_prob();
prob = coco_set(prob, 'tor', 'autonomous', false);
prob = coco_set(prob, 'coll', 'NTST', 40);
prob = coco_set(prob, 'cont', 'NAdapt', 0, 'h_max', 2, 'PtMX', 60);
torargs = \{@vdp @vdp_DFDX @vdp_DFDP @vdp_DFDT t1 x0 {'Om2','c','a','om1','om2','varrho'} 
    [p0' -2*pi/T_po p0(1) -varrho]\};

prob = ode_isol2tor(prob, '', torargs\{:\});

coco(prob, 'vdP_torus', [], 1, {'a','Om2','om2','om1', 'varrho','c'}, \{[0.1 2] [1 2]\});
\end{alltt}
Here we choose negative $\omega_1$ and $\varrho$ because the rotation direction could be in opposite direction. In other words, the value of $\varrho$ can be both positive and negative. One may predict that $\varrho$ is not changed in the continuation run \texttt{vdP\_torus} because releasing the first four parameters $(a,\Omega_2,\omega_2,\omega_1)$ is enough to yield a one-dimensional solution manifold of tori. Indeed $\varrho$ does not change in this continuation run.

Next we use constructor \texttt{ode\_tor2tor} to perform continuation from a saved torus solution in previous run \texttt{vdP\_torus}.
\begin{alltt}
bd  = coco_bd_read('vdP_torus');
lab = coco_bd_labs(bd, 'EP');
lab = max(lab);
prob = coco_prob();
prob = coco_set(prob, 'cont', 'NAdapt', 5, 'h_max', 2, 'PtMX', 60);
prob = ode_tor2tor(prob, '', 'vdP_torus', lab);

coco(prob, 'vdP_torus_varrho', [], 1,
{'a','Om2','om2','varrho','om1','c'}, \{[0.1 2] [1 2]\});
\end{alltt}
In the above continuation run \texttt{vdP\_torus\_varrho}, $\varrho$ is free to change, and several branch points are detected. We then proceed to switch to secondary solution branch passing through one branch point. We use constructor \texttt{ode\_BP2tor} to perform such switch.
\begin{alltt}
bd  = coco_bd_read('vdP_torus_varrho');
lab = coco_bd_labs(bd, 'BP');
lab = lab(1);
prob = coco_prob();
prob = coco_set(prob, 'cont', 'NAdapt', 0, 'h_max', 2, 'PtMX', 60);
prob = ode_BP2tor(prob, '', 'vdP_torus_varrho', lab);

coco(prob, 'vdP_torus_varrho_BP', [], 1,
{'a','Om2','om2','varrho','om1','c'}, {[0.1 2] [1 2]});
\end{alltt}

As exercises to this example, the reader may follow the methodology of previous example. Specifically, continuation of periodic orbits is performed first and then switch to the continuation of tori by \texttt{ode\_TR2tor} and detected TR points in the continuation run of periodic orbits. In addition, the reader may consider the case of positive $\varrho$ and compare the results in two cases.

\section{Discussion}
\label{sec:discussion}
A unified toolbox \texttt{tor} is presented in this report. There are some limitations of the discretization schemes used in the toolbox. First, the number of segments has to be determined before continuation and cannot be changed during continuation. The user needs to check the fidelity of obtained results and then update the number of segments if necessary. An error estimator is instructive for choosing reasonable number of segments. Second, the collocation method used here results in heavy computation load for large dynamical systems. As an alternative, parallel shooting/forward simulation may be utilized to remove such a bottleneck. Finally, the stability of torus solution is not available in current toolbox. Variational equation of torus may be derived and used to determine the stability of torus solution.

\bibliographystyle{plain}
\bibliography{reference}

\end{document}